\documentclass{amsart}

\usepackage{amsmath}
\usepackage{amsthm}
\usepackage{amsfonts}
\usepackage{amssymb}
\usepackage{latexsym}
\usepackage{amscd}
\usepackage{stmaryrd}

\theoremstyle{plain}

\allowdisplaybreaks

\begin{document}

\title{Another proof of Zagier's evaluation formula of the multiple zeta values $\zeta(2,\ldots,2,3,2,\ldots,2)$}

\date{\today\thanks{This work was partially supported by the National Natural
Science Foundation of China (grant no. 11001201) and the Fundamental Research Funds for the Central
Universities. The author is grateful to the referee for his/her useful comments.}}

\author{Zhong-Hua Li}

\address{Department of Mathematics, Tongji University, No. 1239 Siping Road,
Shanghai 200092, China}

\email{zhonghua\_li@tongji.edu.cn}



\maketitle

\begin{abstract}
Using some transformation formulas of the generalized hypergeometric series $\,_3F_2$, we give another proof of
D. Zagier's evaluation formula of the multiple zeta values $\zeta(2,\ldots,2,3,2,\ldots,2)$.
\end{abstract}

\vskip20pt


In a recent paper \cite{zagier}, D. Zagier found and proved an evaluation formula
\begin{align}
\zeta(\underbrace{2,\ldots,2}_a,3,\underbrace{2,\ldots,2}_b)=2\sum\limits_{r=1}^{a+b+1}(-1)^r
c_{a,b}^r\zeta(2r+1)
\zeta(\underbrace{2,\ldots,2}_{a+b+1-r}),
\label{Eq:Zagier}
\end{align}
where $a,b$ are nonnegative integers and
$$c_{a,b}^r=\binom{2r}{2a+2}-\left(1-\frac{1}{2^{2r}}\right)\binom{2r}{2b+1}.$$
Here for positive integers $n,k_1,\ldots,k_n$ with $k_n\geqslant 2$, the multiple zeta value
$\zeta(k_1,\ldots,k_n)$ is defined by the following infinite series
$$\zeta(k_1,\ldots,k_n)=\sum\limits_{0<m_1<\cdots<m_n}\frac{1}{m_1^{k_1}\cdots m_n^{k_n}}.$$
The evaluation formula \eqref{Eq:Zagier} plays an important role in the work of F. Brown \cite{brown},
who showed that all multiple zeta values can be represented as $\mathbb{Q}$-linear combinations of multiple zeta values of the same weight with all arguments are $2$'s and $3$'s and all periods of mixed Tate motives over $\mathbb{Z}$ are
$\mathbb{Q}[(2\pi i)^{\pm 1}]$-linear combinations of multiple zeta values.

As in \cite{zagier}, let $H(a,b)$ (resp. $\widehat{H}(a,b)$) denote the left-hand side (resp. the right-hand side) of \eqref{Eq:Zagier}.
One considers the following two generating functions
\begin{align*}
F(x,y)=&\sum\limits_{a,b=0}^\infty (-1)^{a+b+1}H(a,b)x^{2a+2}y^{2b+1},\\
\widehat{F}(x,y)=&\sum\limits_{a,b=0}^\infty (-1)^{a+b+1}\widehat{H}(a,b)x^{2a+2}y^{2b+1}.
\end{align*}
D. Zagier expressed these two functions by classical special functions. For the function $F$, by \cite[Proposition 1]{zagier}, we know that
\begin{align}
\frac{\pi}{\sin\pi y}\cdot F(x,y)=\left.\frac{d}{dz}\right|_{z=0}\,_3F_2\left({x,-x,z\atop 1+y,1-y};1\right).
\label{Eq:left-function}
\end{align}
On the other hand, by (the proof of) \cite[Proposition 2]{zagier}, we have
\begin{align}
&\frac{\pi}{\sin\pi y}\cdot\widehat{F}(x,y)
\label{Eq:right-function}\\
=&\psi(1+y)+\psi(1-y)-\frac{1}{2}[\psi(1+x+y)+\psi(1-x-y)+\psi(1+x-y)
\nonumber\\
&+\psi(1-x+y)]-\frac{\sin\pi x}{2\sin\pi y}\cdot[\psi(1+(x+y)/2)+\psi(1-(x+y)/2)
\nonumber\\
&-\psi(1+(x-y)/2)-\psi(1-(x-y)/2)-\psi(1+x+y)\nonumber\\
&-\psi(1-x-y)+\psi(1+x-y)+\psi(1-x+y)].
\nonumber
\end{align}
Here the generalized hypergeometric series $\,_3F_2$ is defined as (see \cite{bailey})
$$\,_3F_2\left({\alpha_1,\alpha_2,\alpha_3\atop \beta_1,\beta_2};z\right)=\sum\limits_{n=0}^\infty
\frac{(\alpha_1)_n(\alpha_2)_n(\alpha_3)_n}{n!(\beta_1)_n(\beta_2)_n}z^n,$$
with the ascending Pochhammer symbol
$$(\alpha)_n=\begin{cases}
1, & \text{\,if\;} n=0,\\
\alpha(\alpha+1)\cdots(\alpha+n-1), & \text{\,if\;} n>0.
\end{cases}$$
And $\psi(z)=\Gamma'(z)/\Gamma(z)$ is the digamma function.

D. Zagier proved indirectly that $F=\widehat{F}$.
The purpose of this short note is to give a direct proof of the equality of the right-hand sides of \eqref{Eq:left-function} and \eqref{Eq:right-function}. Our proof uses some transformation
formulas of the $\,_3F_2$-series. To save space, below we will denote the special value $\,_3F_2\left({\alpha_1,\alpha_2,\alpha_3\atop\beta_1,\beta_2};1\right)$ by  $\,_3F_2\left({\alpha_1,\alpha_2,\alpha_3\atop\beta_1,\beta_2}\right)$.
As in \cite{li}, we need two transformation formulas. The first one is (see \cite[Sec. 3.8, Eq. (1), p. 21]{bailey})
\begin{align}
\,_3F_2\left({\alpha_1,\alpha_2,\alpha_3\atop \beta_1,\beta_2}\right)
=&\frac{\Gamma(\beta_1)\Gamma(\beta_1-\alpha_1-\alpha_2)}{\Gamma(\beta_1-\alpha_1)\Gamma(\beta_1-\alpha_2)}\,_3F_2\left({\alpha_1,\alpha_2,\beta_2-\alpha_3\atop \alpha_1+\alpha_2-\beta_1+1,\beta_2}\right)
\label{Eq:trans-to-two}\\
&+\frac{\Gamma(\beta_1)\Gamma(\beta_2)\Gamma(\alpha_1+\alpha_2-\beta_1)
\Gamma(\beta_1+\beta_2-\alpha_1-\alpha_2-\alpha_3)}{\Gamma(\alpha_1)\Gamma(\alpha_2)\Gamma(\beta_2-\alpha_3)
\Gamma(\beta_1+\beta_2-\alpha_1-\alpha_2)}
\nonumber\\
&\times \,_3F_2\left({\beta_1-\alpha_1,\beta_1-\alpha_2,\beta_1+\beta_2-\alpha_1-\alpha_2-\alpha_3\atop \beta_1-\alpha_1-\alpha_2+1,\beta_1+\beta_2-\alpha_1-\alpha_2}\right),
\nonumber
\end{align}
provided that $\Re(\beta_1+\beta_2-\alpha_1-\alpha_2-\alpha_3)>0$ and $\Re(\alpha_3-\beta_1+1)>0$. The second one is (see \cite[Ex. 7, p. 98]{bailey})
\begin{align}
\,_3F_2\left({\alpha_1,\alpha_2,\alpha_3\atop \beta_1,\beta_2}\right)=\frac{\Gamma(\beta_2)\Gamma(\beta_1+\beta_2-\alpha_1-\alpha_2-\alpha_3)}{\Gamma(\beta_2-\alpha_3)\Gamma(\beta_1+\beta_2-\alpha_1-\alpha_2)}\,_3F_2
\left({\beta_1-\alpha_1,\beta_1-\alpha_2,\alpha_3\atop \beta_1,\beta_1+\beta_2-\alpha_1-\alpha_2}\right),
\label{Eq:trans-to-one}
\end{align}
provided that $\Re(\beta_1+\beta_2-\alpha_1-\alpha_2-\alpha_3)>0$ and $\Re(\beta_2-\alpha_3)>0$.

Since
$$(x)_n(-x)_n=\frac{1}{2}(x)_n(1-x)_n+\frac{1}{2}(1+x)_n(-x)_n,$$
we have
\begin{align}
\,_3F_2\left({x,-x,z\atop 1+y,1-y}\right)=\frac{1}{2}\,_3F_2\left({x,1-x,z\atop 1+y,1-y}\right)+\frac{1}{2}\,_3F_2\left({1+x,-x,z\atop 1+y,1-y}\right).
\label{Eq:first-step}
\end{align}
Note that the right-hand side of \eqref{Eq:first-step} is symmetric about $x\leftrightarrow -x$. Hence we only need to consider the first $\,_3F_2$-series of the right-hand side of \eqref{Eq:first-step}. Applying the transformation formula \eqref{Eq:trans-to-two} with
$\alpha_1=x,\alpha_2=z,\alpha_3=1-x,\beta_1=1+y,\beta_2=1-y$, we get
\begin{align}
&\,_3F_2\left({x,1-x,z\atop 1+y,1-y}\right)=\frac{\Gamma(1+y)\Gamma(1-x+y-z)}{\Gamma(1-x+y)\Gamma(1+y-z)}
\,_3F_2\left({x,x-y,z\atop x-y+z,1-y}\right)
\label{Eq:second-step}\\
&\quad+\frac{\Gamma(1+y)\Gamma(1-y)\Gamma(x-y+z-1)\Gamma(1-z)}{\Gamma(x)\Gamma(z)\Gamma(x-y)\Gamma(2-x-z)}
\,_3F_2\left({1-x+y,1+y-z,1-z\atop 2-x+y-z,2-x-z}\right).
\nonumber
\end{align}
Applying the transformation formula \eqref{Eq:trans-to-one} to the first $\,_3F_2$-series of the right-hand side of \eqref{Eq:second-step} with
$\alpha_1=x,\alpha_2=x-y,\alpha_3=z,\beta_1=x-y+z,\beta_2=1-y$, we get
\begin{align}
&\,_3F_2\left({x,1-x,z\atop 1+y,1-y}\right)
\nonumber\\
=&\frac{\Gamma(1+y)}{\Gamma(1+y-z)}\frac{\Gamma(1-y)}{\Gamma(1-y-z)}\frac{\Gamma(1-x-y)}{\Gamma(1-x-y+z)}\frac{\Gamma(1-x+y-z)}{\Gamma(1-x+y)}\nonumber\\
&\qquad\qquad\times\,_3F_2\left({-y+z,z,z\atop x-y+z,1-x-y+z}\right)
\nonumber\\
&\quad+\frac{\Gamma(1+y)\Gamma(1-y)}{\Gamma(x)\Gamma(2-x-z)}\frac{\Gamma(x-y-1+z)}{\Gamma(x-y)}\frac{\Gamma(1-z)}{\Gamma(z)}
\,_3F_2\left({1-x+y,1+y-z,1-z\atop 2-x+y-z,2-x-z}\right)
\nonumber\\
\equiv& [1+\psi(1+y)z][1+\psi(1-y)z][1-\psi(1-x-y)z][1-\psi(1-x+y)z][1+0z]
\nonumber\\
&\quad+\frac{y}{1-x}\frac{\sin\pi x}{\sin\pi y}\frac{z}{x-y-1}\,_3F_2\left({1-x+y,1+y,1\atop 2-x+y,2-x}\right)\qquad\pmod{z^2}
\nonumber
\end{align}
as $z\rightarrow 0$. Here we have used the reflection formula
$$\Gamma(z)\Gamma(1-z)=\frac{\pi}{\sin \pi z}.$$

Now using the summation formula
\begin{align}
\,_3F_2\left({1,\alpha,\beta\atop 1+\alpha,2+\alpha-\beta}\right)=&\frac{\alpha(1+\alpha-\beta)}{\beta-1}
[\psi(\alpha)-\psi(2+\alpha-2\beta)
\label{Eq:summation-formula}\\
&\quad-\psi((\alpha+1)/2)+\psi((\alpha+3)/2-\beta)]
\nonumber
\end{align}
provided that $\Re(\alpha-2\beta)>-2$ with $\alpha=1-x+y,\beta=1+y$,
we get
\begin{align}
&\left.\frac{d}{dz}\right|_{z=0}\,_3F_2\left({x,1-x,z\atop 1+y,1-y}\right)
\label{Eq:last-step}\\
=&\psi(1+y)+\psi(1-y)-\psi(1-x+y)-\psi(1-x-y)
\nonumber\\
&-\frac{\sin\pi x}{\sin\pi y}\cdot[\psi(1-x+y)-\psi(1-x-y)-\psi(1-(x-y)/2)+\psi(1-(x+y)/2)].
\nonumber
\end{align}
Combining \eqref{Eq:first-step} and \eqref{Eq:last-step} gives the desired equality between \eqref{Eq:left-function} and \eqref{Eq:right-function}.

For a proof of the summation formula \eqref{Eq:summation-formula}, one may apply \cite[Eq. (4)]{LGRA} with $m=1$ and
$k=l=0$ to get
\begin{align*}
\,_3F_2\left({1,\alpha,\beta\atop 1+\alpha,2+\alpha-\beta}\right)=&\frac{-\alpha(1+\alpha-\beta)}{2(\beta-1)}
[\psi((\alpha+1)/2)-\psi(\alpha/2)
\\
&\quad+\psi(\alpha/2-\beta+1)-\psi((\alpha+1)/2-\beta+1)].
\end{align*}
The above equation with the help of the formula
$$\psi(2z)=\frac{1}{2}\psi(z)+\frac{1}{2}\psi(z+1/2)+\log 2$$
implies \eqref{Eq:summation-formula}. This completes the proof.



\begin{thebibliography}{99}


\bibitem{bailey} W. N. Bailey, Generalized hypergeometric series, Cambridge University Press, Cambridge, 1935.

\bibitem{brown} F. Brown, Mixed Tate motives over $\mathbb{Z}$, Ann. of Math. 175(2012), 949-976.

\bibitem{LGRA} J. L. Lavoie, F. Grondin, A. K. Rathie and K. Arora, Generalizations of Dixon's theorem on the sum of a $\,_3F_2$,
Math. Comp. 62(1994), 267-272.

\bibitem{li} Z. Li, On a conjecture of Kaneko and Ohno, Pacific J. Math. 257(2)(2012), 419-430.

\bibitem{zagier} D. Zagier, Evaluation of the multiple zeta values $\zeta(2,\ldots,2,3,2,\ldots,2)$, Ann. of Math. 175(2012),
977-1000.

\end{thebibliography}
\end{document}